\documentclass[11pt]{amsart}

\def\factor#1.#2.{\left. \raise 2pt\hbox{$#1$} \right/\hskip -2pt\raise -2pt\hbox{$#2$}}

\def\tbase{\text{Base}\,}

\def\BC{\mathbb C}\def\BF{\mathbb F}

\def\BP{\mathbb P}
\def\pp#1{\mathbb P^{#1}}

\def\pp#1{{\mathbb P}^{#1}}
\def\tdim{\rm dim}
\def\hd{,...,}

\def\cH{{\mathcal H}}

\def\cO{{\mathcal O}}
\def\cP{{\mathcal P}}
\def\CC{\mathbb C}

\def\ZZ{\mathbb Z}
\def\SS{\mathbb S}

\def\11{\mathbf 1}
\def\PP{\mathbb P}
\def\QQ{\mathbb Q}
\def\FF{\mathbb F}

\def\fg{{\mathfrak g}}

\def\l{\lambda}
\def\a{\alpha}

\def\o{\omega}

\def\g{\gamma}
\def\s{\sigma}

\def\m{\mu}

\def\ot{{\mathord{\,\otimes }\,}}
\def\op{{\mathord{\,\oplus }\,}}

\def\lra{{\mathord{\;\longrightarrow\;}}}
\def\ra{{\mathord{\;\rightarrow\;}}}

\def\tdim{\text{dim}\,}

\def\tbase{\text{Base}\,}

\newtheorem{theo}{Theorem}
\newtheorem{coro}[theo]{Corollary}
\newtheorem{lemm}[theo]{Lemma}
\newtheorem{prop}[theo]{Proposition}

\begin{document}

\title{Legendrian varieties}
\author{J.M. Landsberg}
\address{
School of Mathematics,
Georgia Institute of Technology, 
Atlanta, GA 30332-0160, USA}
\author{L. Manivel}
\address{Institut Fourier, UMR 5582 du CNRS, Universit\'e Grenoble I, BP 74,
  38402 Saint Martin d'H\`eres cedex,  FRANCE} 
\begin{abstract}
We investigate the geometry of   Legendrian complex projective 
manifolds $X\subset\PP V$. By definition, this means   $V$ is a complex vector space of dimension 
$2n+2$, endowed with a symplectic form, and the affine tangent space 
to $X$ at each point is a maximal isotropic subspace. We establish basic facts
about their   geometry    and exhibit   examples of inhomogeneous smooth
Legendrian varieties, the first  examples of such in dimension greater than one.
\end{abstract}
\maketitle

\section{Introduction}

The initial motivation for this project stems from the study of the holonomy groups of
Riemannian manifolds, where the only open case for existence of compact
non-homogeneous examples is the quaternion-K\" ahler case. Thanks to work of
Salamon, LeBrun   and others (see, e.g., \cite{sala, LS}),  the question is essentially
equivalent to the existence of inhomogeneous
contact Fano manifolds
(so far none are known). Ye \cite{ye} and others observed that the
set of tangent directions to minimal degree rational curves
through a general point of a contact Fano manifold is a Legendrian
subvariety in its projective span. S. Kebekus  \cite{keb1,keb3} then showed:

\begin{theo}\nonumber  Let $Y$ be a smooth contact Fano manifold with Picard number one, not a 
projective space. Let $y$ a general point of 
$Y$, and denote by $\cH_y\subset\PP T_yY$ the set of tangent directions to 
contact lines on $Y$ passing through $y$. Then $X=\cH_y$ is a smooth
Legendrian variety in its linear span. \end{theo}

Moreover, if  at all points
of $Y$ the corresponding Legendrian variety is homogenous and equivariantly embedded, 
Hong \cite{hong} proved that $Y$ itself must be homogeneous.  

 The homogeneous examples of contact Fano manifolds are as follows:  let $\fg$ denote a  
complex simple Lie algebra and $G$ its adjoint group. Then $G$ has a unique closed orbit
$X_G^{ad}$ inside $\PP\fg$, the projectivization of its Lie algebra -- we call
this variety the {\it adjoint variety} of $G$. The adjoint varieties are contact 
Fano manifolds, and the conjecture of Lebrun and Salamon is that there exists no
other. 

\smallskip
The set of lines passing through a given point of an adjoint variety 
is a smooth homogeneous Legendrian variety.
We call these varieties, in these particular embeddings, the {\it subadjoint varieties}. 
Classical examples are the twisted cubic in $\PP^3$ (coming from the adjoint 
variety of the exceptional group $G_2$), and the products $\PP^1\times\QQ^n$ of 
a projective line with a smooth quadric of dimension $n\ge 1$ (coming from the 
adjoint varieties of the orthogonal groups). Note that the symplectic groups 
give empty subadjoint varieties. The adjoint varieties of the other exceptional 
groups give rise to a  remarkable series of homogeneous varieties, which we called 
the {\it subexceptional series}: they constitute
the third line of the geometric version of Freudenthal's magic square, and were extensively 
studied in \cite{LMfreud, clerc, kaji}; see also \cite{cmr} since they are nice 
examples of varieties with one apparent double point.

\smallskip

In \S 2 we establish a series of Chern class  identities for Legendrian
varieties . They imply that a smooth Legendrian variety which is a product,
must be a $\PP^1\times\QQ^n$ (Corollary \ref{c7}),  that the only Legendrian
embedding of $\PP^n$ is linear (Corollary \ref{c6}),
that a smooth Legendrian variety cannot
be $\pp k$-ruled when $k>1$ (Corollary \ref{cruled}) or
an abelian variety (Corollary \ref{abcor}) and that a homogeneous Legendrian variety
with Picard number one, not necessarily equivariantly embedded a priori, must be a
projective space or a subadjoint variety (Theorem \ref{t11}). But the 
identities  do not exclude, for example, 
that a smooth Legendrian variety has general type, which we do not expect to be possible. 

In \S 3 we establish basic properties about the local differential geometry of
a Legendrian variety $Z$. In particular we show  that any line having contact to
order two with a general point of $Z$ is contained in $Z$ (and thus $Z$ is
uniruled).

\smallskip
Bryant proved that any smooth curve can be embedded in $\PP^3$ as a Legendrian subvariety \cite{bryant}.
 In
contrast, it seems that previous to our work, there were no known examples of
  smooth Legendrian varieties 
of dimension greater than one.

In \S 4 we use Bryant's method to construct   examples of   smooth Legendrian
surfaces, which we now explain in more detail.

\smallskip
By Pfaff's theorem, all (holomorphic) contact structures are locally equivalent
to the space of one-jets of functions and
their Legendrian submanifolds are all locally given by the  one-jets of
functions. In the algebraic category, the model space for the space of
one-jets is $\PP (T^*\PP^{n+1})$ and the Legendrian varieties are
just the lifts $Z^{\#}:=\PP N^*_Z$ (Nash blow-ups) of subvarieties   $Z\subset \pp{n+1}$. Now
$\PP (T^*\PP^{n+1})$ is birational to $\pp{2n+1}$ and one can take a birational
map $\varphi : \PP (T^*\PP^{n+1})\dasharrow \PP^{2n+1}$ that is a linear
contactomorphism on a $\BC^{2n+1}$ (a \lq\lq big cell\rq\rq\ in each space).
 (The inverse rational map and its cousins are studied in detail in \cite{LMclass}.)

 Hence the idea: choose any  variety $Z\subset \PP^{n+1}$, 
then $\tilde{Z}:=\varphi(Z^{\#})$ will be Legendrian in $\PP^{2n+1}$. The problem is that, except for curves, 
this has very little chance to produce a {\it smooth} variety: for example, in dimension two, 
bitangent planes on a surface $Z\subset\PP^3$ tend to produce double points on $\tilde{Z}\subset\PP^5$. 
We analyze   the conditions under which $\tilde{Z}$
can be smooth. 
It turns out that when $Z$ is a {\it Kummer quartic surface} -- a quartic surface with 
sixteen double points as singularities, the resulting surface is smooth. To give a 
precise statement, note that we have a diagram 
$$\begin{array}{rcccl}
 & & Z^{\#} & & \\
 & \swarrow & & \searrow & \\
\PP^3\supset Z & & & & Z^*\subset {\check \PP}^3,
\end{array}$$
where $Z^*\subset {\check \PP}^3$ is the dual variety of $Z$, which is again a Kummer quartic surface 
projectively equivalent to $Z$. The surface $Z^{\#}$ is    a K3 surface, and its natural 
maps to $Z$ and $Z^*$ resolve their singularities. Let $C$ and $D$ denote general hyperplane sections 
of $Z$ and $Z^*$, pulled back to $Z^{\#}$. Then $C$ and $D$ meet transversely in $12$ points. 

\begin{theo} \nonumber The blow up of
the K3 surface $Z^{\#}$  in these twelve points  can be embedded in 
$\PP^5$ as a smooth Legendrian surface of degree $20$. 
\end{theo}

In particular, this provides a counter-example to the na\"\i ve guess that smooth Legendrian 
varieties of dimension greater than one, should be rational.

\section{Chern classes of Legendrian varieties}

In this section we establish Chern class  identities for smooth Legendrian 
varieties. They involve not only the Chern classes of the variety, but also 
the hyperplane class. We determine a number of consequences, including obstructions 
to the existence of a Legendrian embedding of a given variety. For example, 
an abelian variety of dimension at least two has no Legendrian embedding. 

\subsection{An exact sequence}
Let $X\subset \BP V$ be a smooth variety, let 
$x\in X$, let $\hat T_xX\subset V$ denote the affine tangent space 
to $X$ at $x$, let $\tilde T_xX=\BP(\hat T_xX)\subset \BP V$ denote the embedded tangent space and let
$T_xX$ denote the Zariski tangent space.  We have a natural identification 
$$T_xX = Hom(x,\hat T_xX/x)\subset T_x\PP V = Hom(x,V/x).$$
  By hypothesis, $\hat T_xX$ is a Lagrangian 
subspace of $V$, so that the symplectic form induces an identification 
of $V/\hat T_x$ with the dual of $\hat T_x$. 

Consider the  commutative diagram, where $N$ denotes the normal
bundle to $X$:
 
$$\begin{array}{ccccccccc}
 & & & & 0 & & 0 & & \\
 & & & &\uparrow & & \uparrow & & \\
 & & & & \hat T_x^* & = & N_x\ot\cO_x(-1) & & \\
 & & & & \uparrow & & \uparrow & & \\
0 & \rightarrow & \cO_x(-1) & \rightarrow & V\ot\cO_x & \rightarrow
 & T_x\PP V\ot\cO_X(-1) & \rightarrow & 0 \\
 & & || & & \uparrow & & \uparrow & & \\
0 & \rightarrow & \cO_x(-1) & \rightarrow & \hat T_xX & \rightarrow
 & T_xX\ot\cO_x(-1) & \rightarrow & 0 \\
 & & & &\uparrow & & \uparrow & & \\
 & & & & 0 & & 0 & & 
\end{array}$$

We deduce an exact sequence 
$$0\rightarrow \cO_X(-1)  \rightarrow N^*(1) \rightarrow
TX(-1)\rightarrow 0. $$

\subsection{Chern class  identities}
From the previous exact sequence, we get an identity between Chern characters:
$$\begin{array}{rcl}
ch(TX(-1))+ch(\cO_X(-1)) & = & ch(N^*(1)) \\ & =
& ch(\Omega_{\PP V}^1(1)_{|X})-ch(\Omega_X^1(1)) \\
 & = & 2n+2
-ch(\cO_X(1))-ch(\Omega_X^1(1)).
\end{array}$$
Let $h$ denote  the hyperplane class on $X$. We can rewrite the 
previous identity as
$$e^{-h}ch(TX)+e^hch(\Omega_X^1)+e^h+e^{-h}=2n+2.$$
Extracting the homogeneous components, we obtain:

\begin{prop}
Let $X^n\subset\PP^{2n+1}=\PP V$ be a smooth Legendrian variety,
and $h\in H^2(X,\ZZ)$ the hyperplane class. Then for all $m>0$,
the characteristic class
$$\s_{2m}(X,h):=\sum_{i=0}^{2m}(-1)^i\binom{2m}{i}ch_{2m-i}(TX)h^i
+h^{2m}$$
is zero.

In particular, for $m=1$ we get the identity
\begin{equation}\label{cclass1}
2ch_2(TX)=c_1^2-2c_2 = 2hc_1-(n+1)h^2.
\end{equation}
\end{prop}

This already has striking consequences:

\begin{coro}\label{abcor}
An abelian variety, more generally a parallelizable variety, has no Legendrian embedding. 
\end{coro}

\begin{coro}\label{c6}
The unique Legendrian embedding of the projective space $\PP^n$, with $n>1$,
is the linear embedding $\PP^n\subset\PP^{2n+1}$.
\end{coro}

\begin{coro}\label{c7}
Suppose that $X=Y\times Z$ has a Legendrian embedding. 
Then $X=\PP^1\times\QQ^{n-1}$, where $\QQ^{n-1}$ denotes a 
smooth $(n-1)$-dimensional quadric, and the unique Legendrian 
embedding is the Segre embedding.
\end{coro}

\begin{proof} Note that the class $2ch_2=c_1^2-2c_2$ is additive, $ch_2(E\op F)=
ch_2(E)+ch_2(F)$. By the K\"unneth formula, we can decompose
 our very ample class $h$ as $\ell+\theta+m$ over the 
rational numbers, where
$\ell\in H^2(Y,\QQ)$ and $m\in H^2(Z,\QQ)$ are very ample (being the 
classes of the restriction of the hyperplane divisor to fibers of the
projection of $X$ to $Z$ and $Y$, respectively), and 
$\theta\in H^1(Y,\QQ) \ot H^1(Z,\QQ)$. 
Our Chern classes identity decomposes into the following conditions, where  
  $n_Y$ and $n_Z$ respectively denote the dimensions of $Y$ and $Z$:
\begin{eqnarray}\nonumber
c_2(Y) & = & 2\ell c_1(Y)-(n_Y+n_Z+1)\ell^2, \\ \nonumber
0 & = & c_1(Y)\theta - (n_Y+n_Z+1)\ell\theta, \\ \nonumber
0 & = & c_1(Y)m+\ell c_1(Z)- (n_Y+n_Z+1)\ell m, \\ \nonumber
0 & = & \theta c_1(Z)- (n_Y+n_Z+1)\theta m, \\ \nonumber
c_2(Z) & = & 2mc_1(Z)-(n_Y+n_Z+1)m^2.
\end{eqnarray}
Note that the class $\ell'=-c_1(Y)+(n_Y+n_Z+1)\ell$ is very ample 
($\ell$ being very ample, this special case of the 
Fujita conjecture can easily be proved by induction on the dimension). 
By the hard Lefschetz theorem, the
second  identity $\ell'\theta=0$ thus implies that $\theta =0$. 

Now the third condition implies that $c_1(Y)$ (resp. $c_1(Z)$)
and $\ell$ (resp. $m$) are numerically proportional. Let us write
$c_1(Y)=\l \ell$ and $c_1(Z)=\mu m$ for some rational numbers
$\l$ and $\mu$. Then $\l+\m=n_Y+n_Z+1$. Therefore, we cannot
have both $\l\le n_Y$ and $\m\le n_Z$, so we may suppose that 
$\l>n_Y$. Then by the Kobayashi-Ochiai theorem, $Y\simeq \PP^{n_Y}$
is a projective space, $\ell$ is the hyperplane class and $\l=n_Y+1$. 
Hence $\mu=n_Z$, which implies that $Z\simeq\QQ^{n_Z}$ is a 
quadric and $m$ is the hyperplane class. 
But then the first identity cannot be fulfilled, except if $n_Y=1$,
in which case everything vanishes. \end{proof} 

\begin{coro}\label{cruled} Let $X=\BP E$ be the total space of a $\BP^{p}$-bundle over
a variety $Y$. Then if $p>1$, $X$ does not admit a Legendrian embedding.
\end{coro}
\begin{proof}
Let $\pi :\BP E\ra Y$ denote the projection. 
We have exact sequences 
\begin{eqnarray}\nonumber 
0 \lra \cO_X\lra \cO_E(1)\ot\pi^*E^*\lra T^vX\lra 0, \\ \nonumber
0 \lra T^vX\lra TX\lra\pi^*TY\lra 0,
\end{eqnarray}
where $T^vX$ denotes the vertical tangent space with respect to $\pi$.
Let  $\ell$ denote  the first Chern class of the relative hyperplane bundle $\cO_E(1)$, 
and let $r=p+1$,
then $c_1(X)=r\ell+ (basic)$ and $c_2(Z)=\binom r2 \ell^2-\ell.(basic) + (basic)$,\
where $(basic)$ denotes any class that is the pullback   of a class on $Y$.
Suppose that $X$ has a Legendrian embedding, given by a very ample 
line bundle $h=k\ell+\pi^*L$. Let $q=\tdim Y$. Identity \eqref{cclass1} 
implies
$$
r-2rk +(q+r)k^2=0.
$$
Considering this as a quadratic equation for $k$, its discriminant is
$-4rq$.
\end{proof}

\smallskip

\subsection{The case of surfaces}

\begin{prop}
Let $Z\subset\PP^5$ be a ruled Legendrian surface. 
Then $Z=\PP^1\times\PP^1$, embedded by the complete linear 
system $|H+2H'|$. 
\end{prop}

\begin{proof} Suppose that $Z=\PP E$ for some rank $2$ vector bundle $E$ on a curve 
$C$ of genus $g$, and call $\pi$ the projection to $C$. We have exact sequences 
as above. 

Let  $\ell$ denote  the first Chern class of the relative hyperplane bundle $\cO_E(1)$, 
then $c_1(Z)=2\ell-\pi^*(c_1(E)+K_C)$ and $c_2(Z)=\ell^2-\ell.\pi^*(c_1(E)+2K_C)$.
Hence $c_1^2(Z)=2c_2(Z)$.
Note that $\ell^2=deg(E)$, and $\ell.\pi^*L=deg(L)$ for any line bundle $L$ on $C$.

We suppose that $Z$ has a Legendrian embedding, given by a very ample 
line bundle $h=k\ell+\pi^*L$. We must have $2c_1(Z)h=3h^2$, that is, 
$$4k(1-g)=(3k-2)(kdeg(E))+2deg(L)).$$
But $k>0$, and $h^2=k(kdeg(E)+2deg(L))>0$. This implies $g=0$, and then either 
$k=1$ and $deg(E) +deg(L)=1$, or $k=2$ and $deg(E) +2deg(L)=4$. Since $C$ is rational, 
$E$ is split and we can normalize it as $E=\cO\op\cO(-e)$ with $e\ge 0$. Then
the very ampleness of $h$ is equivalent to the condition that $deg(L)>ke$ (\cite{hart}, 
V, Corollary 2.18). 
We easily deduce that $e=0$, so that $Z=\PP^1\times\PP^1$, and then we already 
proved that it must be embedded by the $h=H+2H'$. \end{proof}

\begin{prop}
Let $Z\subset\PP^5$ be a minimal surface of Kodaira dimension zero, in some Legendrian 
embedding. Then   $deg(Z)=8\chi  (\cO_Z)$, and moreover $\chi  (\cO_Z)>1$. 

In particular $Z$ can be neither an abelian nor an Enriques surface. If it is a 
K3 surface, then its genus must be equal to $9$. 
\end{prop}

\begin{proof} The first Chern class of $Z$ is numerically trivial, so that the formula
$c_1^2-2c_2=2c_1h-3h^2$ gives $3deg(Z)=3h^2=2c_2=24\chi  (\cO_Z)$, hence the first claim.
Since $h$ is very ample, Riemann-Roch gives
$$h^0(Z,h)=\chi(h)=\frac{h^2}{2}+\chi  (\cO_Z)=5\chi  (\cO_Z).$$
Since $Z$ is nondegenerate in $\PP^5$, we have $h^0(Z,h)\ge 6$, hence
$\chi  (\cO_Z)>1$. 

As we already know, 
this excludes abelian surfaces, for which $\chi  (\cO_Z)=0$, and also 
Enriques surfaces, 
for which $\chi  (\cO_Z)=1$. For a K3 surface, $\chi  (\cO_Z)=2$, hence $h^2=16$. 
Since the genus of the polarized K3 surface $(Z,h)$ is defined by the identity 
$h^2=2g-2$, we get $g=9$. \end{proof}

\medskip
Following Mukai \cite{mukai2}, a general K3 surface of genus $9$ is a linear section of the symplectic
Grasmmannian $G_{\o}(3,6)$, which we already met as an example of a subadjoint (hence Legendrian)
variety. Note also that a surface in $\PP^5$ is Legendrian, if and only if the image of its Gauss map
is contained in a copy of $G_{\o}(3,6)\subset G(3,6)$. This could   hint  that K3 surfaces of genus
$9$ do admit Legendrian embeddings, but we have not been able to prove or disprove this. 

\begin{prop}
Let $Z\subset\PP^5$ be a minimal surface of general type, in some Legendrian 
embedding. Then $c_1(Z)^2<2c_2(Z)$.
\end{prop}

Note that this is stronger than Miyaoka's bound $c_1(Z)^2\le 3c_2(Z)$. Since Miyaoka's 
bound is sharp, we get infinite families of surfaces of general type admitting 
no Legendrian embedding. 

\medskip\noindent {\it Questions}: Does there exist any example of a surface of 
general type admitting a Legendrian embedding? More generally, what can be the 
Kodaira dimension of a smooth Legendrian variety? Can it be positive? 
\smallskip

\subsection{The case of homogeneous spaces}
 
Let $X=G/P$ be homogeneous, i.e. the quotient of a semi-simple algebraic
group $G$ by some parabolic subgroup $P$. We may suppose that $G$ is 
simple, otherwise $X$ is a product and has already been characterized. 
We fix a maximal torus and a Borel subgroup of $G$ inside $P$, hence a set  $\Delta$ 
of simple roots. Let $I\subset\Delta$ denote  the
simple roots which are not roots of $P$, and $\Phi_P^+$ the 
positive roots which are not roots of $P$. The cohomology algebra 
of $X$ is 
$$H^*(G/P,\ZZ)= \ZZ[\cP]^{W_P}/\ZZ[\cP]^{W}_+\ZZ[\cP]^{W_P},$$
where $\cP$ denotes the root lattice, $W$ the Weyl group of $G$, 
$W_P\subset W$ the Weyl group of $P$, $\ZZ[\cP]^{W_P}$   the 
algebra of $W_P$-invariant polynomials, and  $\ZZ[\cP]^{W}_+\ZZ[\cP]^{W_P}$
the ideal generated by homogeneous $W$-invariants of positive degree 
\cite{hiller}. Since we can interpret  
$\Phi_P^+$ as the set of Chern roots of $TX$, we have 
\begin{eqnarray}\nonumber
c_1(X) = \sum_{\a\in\Phi_P^+}\a, \qquad
2ch_2(X)  = \sum_{\a\in\Phi_P^+}\a^2.
\end{eqnarray}
Since there is no invariant of $W$ in degree one, the identity 
$2ch_2=2c_1l-(n+1)l^2$ means that the corresponding quadratic
element of $\ZZ[\cP]^{W_P}$ is $W$-invariant. But $W$ is generated
by $W_P$ and the simple reflections $s_i$ for $i\in I$, so we just need
to check the invariance under  the action of these simple reflections. 

Recall that a positive root $\a$ is in $\Phi_P^+$ if and only 
$(\a,\o_j)\ge 0$ for all $j\in I$, and there is at least one $k\in I$ 
such that $(\a,\o_k)>0$. For each $i\in I$, define 
$$\Phi^+_P(i) = \{\a\in\Phi^+_P\mid (\a,\o_j)=0 \;  \forall j\in I\backslash i, \;
(\a,\o_i)>0=(s_i\a,\o_i)\}.$$
We have 
$$2ch_2-2s_i(ch_2) = \a_i\sum_{\a\in\Phi_P^+(i)}\a(H_i)(\a+s_i\a).$$

\begin{theo}\label{t11} Let $X=G/P\neq\PP^1$ be a homogeneous space with Picard number one. 
Suppose that $X$ admits a Legendrian embedding, 
not necessarily equivariant a priori. Then $X$ is subadjoint. In particular, the 
Legendrian embedding is the equivariant subadjoint embedding. 
\end{theo}

\begin{proof}  Since 
$X$ has Picard number one, there is a unique simple root $\a_i$ which is
not a root of $P$. Then the corresponding weight $\o_i$ generates $Pic(X)$, 
and we can let $c_1=\g\o_i$ and $\ell=\l\o_i$ for some positive integers
$\g$ and $\l$.  

Consider the identity $2ch_2=2c_1\ell-(n+1)\ell^2$. Applying $s_i$ and dividing
by $\a_i$ we get the relation
$$\sum_{\a\in\Phi_P^+(i)}\a(H_i)(\a+s_i\a)=\l(2\g-(n+1)\l)(2\o_i-\a_i).$$
In particular, the scalar product with $\o_i$ gives 
$$\sum_{\a\in\Phi_P^+(i)}\a(H_i)^2\frac{(\a_i,\a_i)}{2}=\l(2\g-(n+1)\l)(2\o_i-\a_i,\o_i)>0.$$
But $2\o_i-\a_i$ is a linear combination, with positive coefficients, of the fundamental 
weights corresponding to the nodes of the Dynkin diagram that are connected to $i$. 
In particular, $(2\o_i-\a_i,\o_i)>0$ and we deduce 
that $2\g-(n+1)\l>0$. This implies that $\l=1$ (otherwise the index of $X$ would be
greater that $n+1$, which is impossible),  and that $X$ has index $\g>\frac{n+1}{2}$. 

The index of the rational homogeneous spaces with Picard number one 
have been computed by Snow \cite{snow}, 
and we easily get  the following lemma from his results. 

\begin{lemm}
Let $X$ be a homogeneous space with Picard number one. 
\begin{enumerate}
\item If $X$ is adjoint, its index
is $\g=\frac{n+1}{2}$. 
\item If $\g>\frac{n+1}{2}$, then $X$ is either a projective space, 
$G(2,n)$, $G(3,7)$, $G_{\o}(2,n)$,  a quadric, a spinor variety 
$\SS_m$ with $m\le 7$, the Cayley plane, or a subadjoint variety. 
\end{enumerate}\end{lemm}

Then we check that our Chern class identity only holds in the subadjoint case. 
Moreover, $\l=1$ implies that the embedding is given by the subadjoint embedding, 
which is Legendrian, possibly followed by a projection, which excluded by the fact that 
 the dimension of the ambient space would become too small. \end{proof}

\subsection{More Chern class  identities}
If the dimension of $X$ is large enough, we can eliminate the hyperplane class from the
identities $\s_{2m}(X,h)=0$ to obtain identities between the Chern classes of $X$. 
This illustrates the principle that, the greater   the dimension, the more difficult it
should be to find inhomogeneous smooth Legendrian varieties -- if any exist. 

Formally, we just interpret  our identities as polynomial equations for $h$, with coefficients
in the Chow ring of $X$. Let $R_{l,m}$ be the resultant of the polynomials $\s_{2l}$
and  $\s_{2m}$, a polynomial of degree $2l+2m+2$ in the Chern classes. 

\begin{prop} Let $X$ be a smooth variety of dimension $n\ge 8$, admitting a 
Legendrian embedding. Then $R_{l,m}(X)=0$ for $1\le l<m$.
\end{prop}

The first identity $R_{1,2}(X)=0$ comes out in degree $8$ and can be written 
$$C_{8,4}(n+1)^4+C_{8,3}(n+1)^3+C_{8,2}(n+1)^2+C_{8,1}(n+1)+C_{8,0}=0,$$
where the classes $C_{8,i}$ are expressed as follows in terms of the Chern classes:
\begin{eqnarray}\nonumber
C_{8,4} &= & ch_4^2, \\ \nonumber
C_{8,3} &= & 16ch_2ch_3^3-8ch_1ch_3ch_4-20ch_2^2ch_4, \\ \nonumber
C_{8,2} &= & 32ch_1^2ch_2ch_4-16ch_1ch_2^2ch_4+100ch_2^4, \\ \nonumber
C_{8,1} &= & 32ch_1^3ch_2ch_3-176ch_1^2ch_2^3-16ch_1^4ch_4, \\ \nonumber
C_{8,0} &= & 468ch_1^4ch_2^2.
\end{eqnarray}

\section{Local differential geometry of Legendrian varieties}

Let $X^n\subset \BP V$ be a Legendrian variety and let $z\in X$ be a general
point.
By Pfaff's theorem (see, e.g., \cite{BCGGG} p38) there exist local coordinates
$x^1\hd x^{2n+1}$ about $z$ such that locally $X$ may be written as a graph  
of the form
\begin{align} 
 x^{2n+1}&=f(x^1\hd x^n)\\
x^{n+j}&= \frac{\partial f}{\partial x^j}(x^1\hd x^n) \ \ 1\leq j\leq n
\end{align}
where $f$ and all its first derivatives vanish at $z=(0\hd 0)$. 
 Thus $|III_z|$, the third fundamental form, consists of a single cubic
(given by the third order partials of $f$ at $z$)
and the quadrics of the
second fundamental form $|II_z|$ consist  of the $n$ partial derivatives of $|III_z|$.
 In particular, letting $N_{2,z}=II_z(S^2T_zX)$ denote the image of the second fundamental
form, 
we may view $II_z\in  S^2T_z^*\ot N_{2,z}$  as having an additional symmetry,
using the symplectic form $\o$ to identify $T_z$ with $N_{2,z}$ we obtain
$II_z\in  S^3T_z^*$. 

A coordinate free way to see this symmetry is as follows:
Let $\gamma: X\ra G(n+1,2(n+1))$ denote the Gauss map of $X$. Then $X$ is
Legendrian iff the image of $\gamma$ is contained
in the $\o$-isotropic Grassmanian $G_{\o}(n+1,2(n+1))$. Recall that
$T_EG_{\o}(n+1,2(n+1))\simeq S^2E^*$ (see, e.g., \cite{LM0})
so $II_z\simeq d\gamma_z \in T^*_zX\ot T_{T^*_zX}G_{\o}(n+1,2(n+1))
= T^*_zX\ot S^2T^*_zX$. On the other hand, the differential of the
Gauss map of any variety is an element of $S^2 T^*_zX\ot  T^*_zX$ and
the intersection of these two spaces is exactly $S^3T^*_zX$.

This property propagates to higher order differential
invariants. Let $F_{k,z}\in S^kT^*_zX\ot N_{2,z}X$ denote the relative differential
invariant of order $k$ (defined in coordinates in the proof below, see \cite{IvL}  section 3.5  for a 
coordinate free definition):

\begin{prop}\label{ff} Let $X\subset \BP V$ be Legendrian, let $z\in X$ be a smooth
point.  Then using the identification $N_{2,z}X\simeq T^*_zX$,
we have $F_{k,z}\in S^{k+1}T^*_zX$.
\end{prop}

\begin{proof} The coefficients of $F_{k,z}$ are simply the higher order
terms in the Taylor series for the $x^{n+j}$, which correspond to
the $k$-th derivatives of $\frac{\partial f}{\partial x^j}$.
\end{proof}

If $P$ denotes the cubic for the third fundamental form and $v\in T_zX$, let
$Q_v=\frac{\partial P}{\partial v}$ denote the corresponding quadric in $II_z$.
Let $\tbase |II_{z}|:=\{ v\in T_zX\mid II(v,v)=0\}$, the {\it base locus} of the 
second fundamental form, the set of tangent directions to lines having contact
with $X$ at $z$ to order two.

\begin{coro} $\tbase |II_{z}|=\{ v\in T_z \mid v\in (Q_v)_{sing}\}$
\end{coro}

Let $C_z\subset  T_zX$ denote the tangent directions to lines on $X$ passing
through $z$. Note that one always has $C_z\subseteq \tbase |II_{z}|$.

\begin{theo} Let $X\subset \BP V$ be Legendrian and let $z\in X$ be
a general point. Then $C_z=\tbase |II_{z}|$. In particular,  
$ |III_{z}| $ is singular if and only if $X$ is uniruled by lines.
\end{theo}

\begin{proof} Let $v\in \tbase |II_z|$, so $v\in (Q_v)_{sing}$. Now
\cite{Lhss} (3.1.2) implies $v\in (F_{3,v})_{sing}$ (here we
mean the cubic in the normal direction corresponding to the tangent
vector $v$) but now the symmetry
implies $v\in \tbase F_3$. But now \cite{Lhss} (3.1.3) implies $v\in (F_{4,v})_{sing}$
and the symmetry again implies $v\in \tbase F_4$. Continuing,  
$v\in (F_{j,v})_{sing}$ and $v\in \tbase (F_i)$ for all $i\leq j$ implies
$v\in (F_{j+1,v})_{sing}$ and then the symmetry implies in turn 
$v\in \tbase (F_{j+1})$, and these two facts imply $v\in (F_{j+2,v})_{sing}$
etc... and one obtains $v\in \tbase(F_l)$ for all $l$, i.e., that there
is a line having infinite order contact to $X$ at $x$ in the direction of $v$.
\end{proof}

We summarize some other known properties of Legendrian varieties:
 
\begin{prop} Let $X\subset\PP^{2n+1}=\PP V$ be a   Legendrian variety.
\begin{enumerate}
\item If $X$ is linearly degenerate, then $X$ is a linear subspace or a cone 
over a linear subspace. 
\item If $X$ is not linearly degenerate, then the tangent variety $\tau(X)\subset\PP V$ 
and the dual variety $X^* \subset\PP V^*$ are projectively isomorphic hypersurfaces.
\end{enumerate}
\end{prop}

\begin{proof} 
 If $X$ is degenerate, any tangent space to $X$ is contained in some 
 hyperplane. By duality, we get that every tangent space to $X$ passes through 
 some fixed point. In characteristic zero, this implies that $X$ is a cone.

A hyperplane $H$ is tangent to $X$ at the point $x$
if and only if the point $h=H^{\perp}$ belongs to the embedded tangent space $T_xX$. 
This gives a linear identification between $X^* $ and $\tau(X)$. 
The fact that they are   hypersurfaces  follows e.g.,
from the fact that $III_z$ is nonzero (because
$X$ is not linearly degenerate) hence, by the infinitesimal calculation of
the dimension of the secant variety $\s(X)$ in \cite{GH2}, we have 
$\s (X)=\BP V$ and thus by the Fulton-Hansen connectedness theorem $\tau(X)$ is a hypersurface. 
\end{proof} 

\medskip In the case of surfaces, the degree of the dual variety, which is sometimes 
called the {\it codegree}, does not depend of the embedding but only of the Chern numbers. 
This seems to be specific to dimension two. Note the curious relation with Miyaoka's 
inequality. 

\begin{prop} If $X\subset\PP^5$ is a Legendrian surface, its codegree is equal to  
$3c_2-c_1^2$. \end{prop}

\begin{proof} By Katz's formula (\cite{GKZ}, Chapter 2, Theorem 3.4), 
the degree of $X^*$ is $c_2-2hc_1+3h^3$. Since for
a Legendrian surface we know that $c_1^2-2c_2=2hc_1-3h^3$, the result follows. \end{proof}

\section{Bryant's method}

Now that we have various numerical conditions on smooth Legendrian varieties,
we expect that it should be a   delicate problem to construct explicit examples,
especially in higher dimensions. The only method of construction that we are aware of
was suggested by Bryant. As explained in the introduction, it is based on the 
observation that the corresponding problem in $\PP(T^*\PP^{n+1})$ is easily solved. 
One then tries to transport the solutions to $\PP^{2n+1}$ through an explicit birational
map. Unfortunately, there are  strong obstructions for this simple and elegant idea
to produce smooth varieties and we are only able to make it work in a very special
situation.

\subsection{The birational map}
For symmetry reasons we consider $\PP(T^*\PP^{n+1})$ as the  flag variety 
$\FF_{1,n}(\CC^{n+1})$ of pairs of incident lines and hyperplanes in $\PP^n$. It 
has two projections $p$ and $p'$ on $\PP^n$ and its dual $\check{\PP}^n$. 
If we choose homogeneous coordinates $x_0,\ldots ,x_n$ on $\PP^n$, and
dual coordinates  
$y^0,\ldots ,y^n$    on $\check{\PP}^n$, it is 
the subvariety of $\PP^n\times \check{\PP}^n$ defined by the equation 
$\sum_{i=0}^nx_iy^i=0$.

\smallskip
Let $[w_1,\ldots,w_n,z_1,\ldots ,z_n]$ be homogeneous 
coordinates on $\PP^{2n-1}$.
We consider the rational map 
$$\begin{array}{rrcl}
\varphi : & \FF_{1,n}(\CC^{n+1}) & \dasharrow & \PP^{2n-1} \\
          &   ([x],[y]) & \mapsto & [x_0y^1,\ldots,x_0y^{n-1},x_0y^0-x_ny^n,x_{n-1}y^n,\ldots , x_0y^n,
x_1y^n].\end{array}$$
This is a birational map, whose inverse is given by
$$\begin{array}{rrcl}
\varphi^{-1} : & \PP^{2n-1} & \dasharrow & \FF_{1,n}(\CC^{n+1}) \\
 &  [w,z]& \mapsto & ([z_n,z_1,\ldots ,z_{n-1},-\frac{1}{2}(w_n+\frac{(w,z)}{z_n})], \\
 & & & \hspace{3cm} [\frac{1}{2}(w_n-\frac{(w,z)}{z_n}), w_1,\ldots ,w_{n-1},z_n]),
\end{array}$$
where $(w,z):=\sum_{i=1}^{n-1} w_iz_i$.

\smallskip
Consider the contact structures given by the line-bundle valued
one-forms $\theta' = xdy=-ydx$ on $\FF_{1,n}(\CC^{n+1})$, and 
$\theta=zdw-wdz$ on $\PP^{2n-1}$.

\begin{lemm} The birational map $\varphi$ is compatible with the contact structures
determined by $\theta '$
on $\FF_{1,n}(\CC^{n+1})$ and $\theta$ on $\PP^{2n-1}$.
\end{lemm}

\begin{proof} A simple computation shows that $\varphi^*\theta=x_0y^n\theta'$. \end{proof}

\begin{coro} 
Let $Z\subset\PP^n$ be a subvariety, and $Z^{\#}\subset\FF_{1,n}(\CC^{n+1})$ its conormal
variety. Then $\tilde{Z}=\varphi(Z^{\#})$ is a Legendrian subvariety of $\PP^{2n-1}$.
\end{coro}

Recall that the conormal variety is defined as the closure of
the incidence variety of pairs $(z,H)$ where $z\in Z$ is a smooth 
point and $H\in \pp{n*}$ is such that $\tilde T_zZ\subset H$.   Its projection to the dual 
projective space ${\check\PP}^n$ is, by definition, the dual variety $Z^*$ of $Z$. 
That the image variety $\tilde{Z}$ is Legendrian means that at every smooth point, the
affine tangent space is maximal isotropic with respect to   $\theta$. 
Unfortunately, the fact that $\varphi$ is not an isomorphism will tend to produce 
singularities on $\tilde{Z}$. We now analyze $\varphi$ in some detail to determine
conditions under which $\tilde Z$ is smooth.

\bigskip\noindent {\bf Fact 0}. {  The exceptional locus $Exc(\varphi)$
of $\varphi$ is the hyperplane section $x_0y^n=0$ of $\FF =\FF_{1,n}(\CC^{n+1})$, union of the two 
irreducible divisors $E_1:=\{ ([x],[y])\mid x_0=0\}$ and $E_2:=\{ ([x],[y])\mid y^n=0\}$. 
The indeterminacy locus is $Ind(\varphi)=E_1\cap E_2$. 

Outside the exceptional locus, $\varphi$ restricts to an isomorphism
onto the affine space defined as the complement of the hyperplane 
$P=(z_n=0)$.}

\bigskip Let $H_0$ be the hyperplane $\{ x_0=0\}$ in $\PP^n$, and 
$p_0\in H_0$ the point dual to the hyperplane $y^n=0$ in ${\check\PP}^{n}$.
Geometrically,   $E_2=\{ (p, H)\in \FF \mid p_0\in H\}$
and  $E_1= \{ (p, H)\in \FF \mid p \in H_0\}$. 

\bigskip\noindent {\bf Fact 1}. {  Two points $(p, H)$ and $(p', H')$ of $\FF$
outside the indeterminacy locus are in the same fiber of $\varphi$
if and only if}
$$p=p'\in H_0 \quad or \quad p_0\in H=H'.$$

\begin{proof} On $Exc(\varphi)-Ind(\varphi)$, $\varphi$ is given by the 
following formulas:
$$\begin{array}{rcl}
\varphi([0,x_1,\ldots ,x_n],[y]) & = & [0,\ldots ,0,x_n,\ldots ,x_1,0], \\
\varphi([x],[y^0,\ldots ,y^{n-1},0]) & = & 
[y^1,\ldots ,y^{n-1},-y_0,0,\ldots ,0].
\end{array}$$ 
This implies the claim. \end{proof}

\bigskip Let $Z\subset\PP^n$ be some irreducible, possibly singular
hypersurface. 

\bigskip\noindent {\bf Fact 2}. {   If $\tilde Z$ is smooth, then
$(Z-Z\cap H_0)^{\#}$ is smooth
outside $E_2$.}

\begin{proof} This follows from Fact 0, since $$Z^{\#}-Exc(\varphi)=
(Z-Z\cap H_0)^{\#}-E_2.$$\end{proof}

\medskip We say that $Z$ is {\sl quasi-smooth} outside $H_0$ if 
$(Z-Z\cap H_0)^{\#}$ is smooth.  For a curve, this means that outside $H_0$, 
$Z$ has only nodes or simple cusps. In general, its seems to be a difficult 
problem to understand quasi-smooth singularities. 
A surface with a double curve will be quasi-smooth 
at smooth points of that curve. An isolated quadratic singularity
is also quasi-smooth, but apparently no other simple surface singularity. 
Arnold and its school have classified what they call (real) stable 
Legendrian singularities (\cite{arnold}, section 21). 
This gives examples of quasi-smooth singularities, for example
the swallow-tail. But note that $Z$ and its dual have the same conormal
variety, so that   the dual variety of a smooth variety
(which is in general very singular), is always quasi-smooth.

\bigskip\noindent {\bf Fact 3}. {Let $z\in Z\cap H_0$ be a smooth 
point such that $p_0\notin T_zZ$. Then $Z^{\#}$ meets  the fiber 
of $\varphi$ containing $(z,\tilde T_zZ)$ transversely at $(z,\tilde T_zZ)$.}
 
\begin{proof}
By Fact 1, the fiber of $\varphi$ at $(z,\tilde T_zZ)$ is
$\{ (z,H)\mid z\in H\}$, or, identifying  a linear space on $\BF$ 
with a vector subspace of the ambient $V\ot V^*$, it is the set of elements
$\{ u\ot h' \mid [u]=z, \langle u,h'\rangle =0\}$, which also corresponds
to the kernel of $\varphi_*$ at $(z,\tilde T_zZ)$.  On the other hand, 
an easy moving frames calculation shows
$$
\hat  T_{z,\tilde T_zZ}Z^{\#}  =
\{ v\ot h - x\ot   Q(v,\cdot)   \mid [h]=\hat T_zZ,\  v\in \hat T_zZ\} \subset V\ot V^*  
$$ 
where $Q\in S^2T^*$ is the quadric corresponding to $h$ as an
element of the second fundamental form
of $Z$ at $z$. (Here we slightly abuse notation to consider
$Q(v,\cdot)$ as an element of $V^*$ which requires a choice of splitting,
but this identification is harmless.)
The unique intersection of these two linear spaces is the line $\{ u\ot h
\mid [u]=z,\ [h]=\tilde T_zZ\}$, which corresponds to the zero vector in $T_{z,\tilde T_zZ}Z^{\#}$.
\end{proof}

\bigskip\noindent {\bf Fact 3'}. {Let $z\in Z-Z\cap H_0$ be a smooth 
point such that $p_0\in T_zZ$. Then $Z^{\#}$ meets  the fiber 
of $\varphi$ containing $(z,  T_zZ)$ transversely at $(z,\tilde T_zZ)$, if and only if the Gauss 
map of $Z$ is immersive at $z$, that is, if and only if $z$ is not a 
flex point of $Z$.}

\begin{proof} Here, by Fact 1 again, the  fiber of $\varphi$ containing $(z,\tilde T_zZ)$
is the projective space $\{(p,\tilde T_zZ)\mid p\in \tilde T_zZ\}$ which    we
identify with the linear subspace  $\{ v\ot h \mid [h]=\tilde T_zZ,\ \langle v,h\rangle =0\}$
of
$V\ot V^*$, which also corresponds
to the kernel of $\varphi_*$ at $(z,\tilde T_zZ)$.
Note 
in particular that this space consists of rank one elements. On the other hand, as above
$$
\hat  T_{z,T_zZ}Z^{\#}  =
\{ v\ot h - x\ot   Q(v,\cdot)  \mid [h]=\hat T_zZ,\ v\in \hat T_zZ\}   
$$ 
which, except for the point
$v\ot h$ when $[v]=z$  consists of rank at least two elements, as long as
  $v$ is not a singular point of $Q$. If $v$
is a singular point of $Q$ then the two spaces coincide.
\end{proof}

\bigskip\noindent {\bf Fact 3''}. { More generally, still assuming $\tilde Z$ is smooth, let 
$z\in Z-Z\cap H_0$, not necessarily a smooth point,  and $z^*$ a tangent
hyperplane at $z$, containing $p_0$. 
Suppose  that $Z^{\#}$ is smooth at $(z,z^*)$. 

Then $Z^{\#}$ meets the fiber of $\varphi$  transversely at $(z,z^*)$, 
if and only if the projection $Z^{\#}\lra Z^{\ast}$   is immersive at $(z,z^*)$. }

\begin{proof} Same proof as for Fact 3'.\end{proof}

\smallskip 
Note that from these facts one can reproduce Bryant's proof that
any smooth projective curve
can be embedded as a smooth Legendrian curve (\cite{bryant}, Theorem G). Simply project the curve
to $\pp 2$ so it has only nodal singularities and make sure it is  in
sufficiently general position with respect to $(p_0,H_0)$. See \cite{bryant} for
details.

\bigskip If $Z$ has dimension greater than one, $Z\cap H_0$ is a positive
dimensional hypersurface in $H_0$. Its tangent hyperplanes 
will   cover $H_0$ (as long as $Z\cap H_0$ is not set theoretically
a linear space), and in particular some of them will
contain $p_0$, so that $Z^{\#}$ will meet $Ind(\varphi)$. 
Before resolving the indeterminacies of $\varphi$, we make two 
simple observations.  

\bigskip\noindent {\bf Fact 4}. { Assume $\tilde Z$ is smooth. Let $p,p'$ be smooth points of
$Z-Z\cap H_0$ with the same tangent hyperplane $H$, passing through $p_0$. 
Then $H$ must be tangent to $Z$ along a curve.} 

\begin{proof} Otherwise $(Z-Z \cap H_0)^{\#}$ would meet a fiber of $\varphi$ 
along a disconnected subset, and its image would therefore be
multibranch at the corresponding point of $\PP^{2n-1}$, 
contradicting the smoothness assumption. \end{proof}

\bigskip\noindent {\bf Fact 4'}. { Assume $\tilde Z$ is smooth. Let $p\in H_0$ be a singular 
point of $Z$. Then $Z$ has at most one branch at $p$ tangent to a 
hyperplane not containing $p_0$.}

\begin{proof} Otherwise $Z^{\#}$ would contain two points $(p,H)$ and $(p,H')$
in the same fiber of $\varphi$ (or if $H=H'$, $Z^{\#}$ would not be 
unibranch at $(p,H)$), and its image under $\varphi$ would not be unibranch, 
thus would be singular. \end{proof} 

\bigskip\noindent To go further in our analysis, we need to resolve the indeterminacies
of $\varphi$. A simple thing to do would be to blow-up the indeterminacy locus, but:

\bigskip\noindent {\bf Fact 5}. {  The indeterminacy locus 
$Ind(\varphi)$ has a quadratic singularity at $(p_0, H_0)$, 
and is smooth outside that point.}

\begin{proof} We choose local coordinates on $\FF:=\FF_{1,n}(\CC^{n+1})$ at $(p_0, H_0)$ as 
follows: for a pair $(p, H)$, $p=[x_0,x_1,\ldots ,x_{n-1},1]$ 
and $H$ is the hyperplane generated by $p$ and $n-1$ other vectors
$e_{n-1}-z_{n-1}e_0,\ldots , e_1-z_1e_0$, so that 
$H=[1,z_1,\ldots ,z_{n-1},-x_0-(x,z)]$, where $(x,z)=\sum_{1\le i\le n-1}x_iz_i$. 
Then the condition that 
$p\in H_0$ is equivalent to $x_0=0$, and the condition that 
$p_0\in H$ is equivalent to $x_0+(x,z)=0$. We thus see that 
$E_1$ and $E_2$ are smooth hypersurfaces meeting nontransversely 
at $(p_0, H_0)$ along the codimension two subvariety $x_0=(x,z)=0$,
a quadratic cone in a coordinate hyperplane. \end{proof}

\bigskip Outside this singularity there is no serious problem:
blowing-up the indeterminacy locus is enough to resolve the
indeterminacies. 

\bigskip\noindent {\bf Fact 6}. {Let $\sigma : \tilde{\FF}
 \ra\FF-\{p_0\in H_0\}$ be the blow-up of $Ind(\varphi)-\{p_0, H_0\}$. 
Then $\psi:=\varphi\circ\sigma$ is a morphism.}

\begin{proof} We check this in local coordinates. Let $(p_1, H_1)$ be a 
point of $Ind(\varphi)$. We may suppose that $p_1=[0,1,0,\ldots ,0]$
and $H_1=[0,\ldots ,0,1,0]$. If $[x_0,1,x_2,\ldots ,x_n]$ are affine
coordinates around $p_1$ in $\PP^n$, and $[y^0,\ldots ,y^{n-2},1,y^n]$
are affine coordinates around $H_1$ in the dual projective space,
we can choose $x_0,x_2, \ldots ,x_n,y^0,y^2,\ldots ,y^{n-2},y^n$
as affine coordinates on $\FF$,  with the missing coordinate $y^1$   
given by the relation $$x_0y^0+y^1+\sum_{i=2}^{n-2}x_iy^i+x_{n-1}
+x_ny^n=0.$$ Above the corresponding open set $U\subset\FF$, the 
blow-up $\tilde{\FF}$ is the set of points $(x,y,[s,t])\in U\times\PP^1$ 
such that $sx_0=ty^n$.  

The roles of $s$ and $t$ being symmetric, we may suppose that $s\neq 0$. 
Then we let $s=1$ and choose 
$x_2, \ldots ,x_n,y^0,y^2,\ldots ,y^{n-2},y^n,t$ as local 
coordinates on $\tilde{\FF}$. In these coordinates, $\psi$ is given by
$$\begin{array}{rcl}
\psi(x,y,[s,t]) & = & \varphi([ty^n,1,x_2,\ldots,x_n],
 [y^0,y^1,\ldots ,y^{n-2},1,y^n] \\
& = & [ty^ny^1,\ldots ,ty^ny^{n-1},x_ny^n-ty^ny^0, \\
 & & \hspace{2.6cm} x_{n-1}y^n,
\ldots ,x_2y^n,y^n,t(y^n)^2] \\
 & = & [ty^1,\ldots ,ty^{n-1},x_n-ty^0,x_{n-1},\ldots, x_2,1,ty^n]. 
\end{array}$$
This is always defined, and our claim follows.\end{proof}

\bigskip Although we do not use it in what follows, we briefly  
describe what happens 
near the singular point. We leave the details to the reader. 

We begin by blowing-up $(p_0, H_0)$ in $\FF$,
giving a map $\sigma : \FF'\ra\FF$. We let $\varphi'=\varphi\circ
\sigma :\FF'\dasharrow\PP^{2n-1}$, and $I'$ be the strict 
transform of $Ind(\varphi)$. 

\bigskip\noindent {\bf Fact 7}. {  The indeterminacy locus 
$Ind(\varphi')$ has two irreducible components: $I'$ and a
smooth hyperplane $H$ inside the exceptional divisor $E\simeq 
\PP^{2n-2}$. These two components meet transversely along the
smooth subvariety $I'\cap H$.}

\bigskip Then we blow-up $I'$ by $\sigma ' :\FF''\ra\FF'$, and we let 
$\varphi''=\varphi'\circ\sigma '$. 

\bigskip\noindent {\bf Fact 8}. {  The indeterminacy locus 
$Ind(\varphi'')$ is the strict transform $H'$ of $H$, a smooth
subvariety of $\FF''$ of codimension two.}

\bigskip Finally we blow-up $H'$ by $\sigma '' :\FF'''\ra\FF''$, and we let 
$\varphi'''=\varphi''\circ\sigma ''$. 

\bigskip\noindent {\bf  Fact 9}. {  $\varphi''':\FF'''\ra\PP^{2n-1}$
is a morphism, and its exceptional locus is the union of four 
irreducible divisors.}

\subsection{Surfaces}

\bigskip\noindent Now suppose that $Z$ is a surface in $\PP^3$, such that 
$Z^{\#}$ does not contain $(p_0, H_0)$ -- for example, we can just ask that 
$p_0\notin Z$. Then by Fact 6, we just need to blow-up 
the smooth part of $Ind(\varphi)$ to resolve the indeterminacies of $\varphi$
on $Z^{\#}$. 

Let $Z^{\#\#}$ denote the strict transform of $Z^{\#}$ by the blow-up
$\sigma$. Also we let $E_1^{\#}$ and $E_2^{\#}$ be the strict transforms of
$E_1$ and $E_2$, and $E_0^{\#}$ be the exceptional divisor of $\sigma$. The
morphism $\psi$ restricts to an isomorphism between the complement of 
$E_0^{\#}\cup E_1^{\#}\cup E_2^{\#}$ and the complement of the hyperplane $H$ in 
$\PP^5$. Moreover, $\psi$ maps $E_0^{\#}$ to a quadratic cone $Q$ inside
the hyperplane $P=\{ z_3=0\}$, 
and $E_1^{\#}$, $E_2^{\#}$ to planes $U_1$ and $U_2$ inside $Q$, meeting 
at the vertex of $Q$. 
 
\begin{prop}\label{p17} Suppose that: \begin{enumerate}
 \item $Z\cap H_0$ is a smooth curve, 
 \item $Z$ has no bitangent plane containing $p_0$, 
 \item $Z\cap H_0$ has no bitangent line containing $p_0$.
 \item $Z^{\#}$ meets $Ind(\varphi)$ transversely.
\end{enumerate}
Then $\psi$ is injective on $Z^{\#\#}$.
\end{prop}

Actually, we can replace condition $(3)$ by the weaker condition $(3')$: 
on a line in $H_0$ containing $p_0$, there is at most one point of $Z$ whose
tangent hyperplane contains $p_0$. 

\begin{proof} We just need to check that the fiber of a point of the quadratic cone $Q$ 
contains at most one point of $Z^{\#\#}$. 

Let $q=[u^1,u^2,w,v_2,v_1,0]$ be a point of 
$Q-U_1\cup U_2$. Let $p\in \pi\circ\sigma(\psi^{-1}(q))\subset Ind(\varphi)$. 
By the formulas above for the morphism $\psi$ we see that $p$ belongs
to the line in $H_0$ generated  by $v=[0,v_1,v_2,0]$ and $p_0$. 
By $(1)$, $Z$ is smooth at $p$, and by $(3)$ or $(3')$, $p$ is uniquely 
determined.  
Now, $(4)$ guarantees that the restriction of $\s$ to $Z^{\#}$ is just the 
blow-up of the finite set of points on the smooth curve $C=Z\cup H_0$, 
where the tangent lines to $C$ hits $p_0$. Our point $p$ is one of these, 
and it follows that $Z^{\#\#}$ is smooth over $p$ and contains the whole 
fiber of $\psi$ over $(p, T_pZ)$. This fiber is a projective line. But it 
follows from its expression that $\psi$ restricts to an isomorphism 
between that line and a line in $\PP^5$. We conclude that  $\psi$
restricted to $Z^{\#\#}$ is injective over $Q-U_1\cup U_2$.
 
Because of Fact 1, $(2)$ and $(3)$ ensure that it is also injective 
over $U_1\cup U_2$, and our claim follows. \end{proof}

\begin{prop}\label{p18} Suppose moreover that: 
\begin{enumerate}
 \item $Z$ is quasi-smooth,
 \item the projection $Z^{\#}\lra Z^*$ is immersive over any 
tangent plane containing $p_0$, 
 \item at any point $p\in Z\cap H_0$ such that $T_pZ$ contains
$p_0$, the base locus of the second fundamental form of $Z$ does not 
contain the line $\overline{pp_0}$.
\end{enumerate}
Then $\psi$ defines an embedding of $Z^{\#\#}$ in $\PP^5$, 
and its image $\tilde{Z}$ is a smooth Legendrian surface in $\PP^5$. 
\end{prop}

\begin{proof} Since $Z^{\#}$ is smooth and meets $Ind(\varphi)$ 
transversely, $Z^{\#\#}$ is smooth (note that this transversality
hypothesis is equivalent to the hypothesis we made on the second 
fundamental form, as the computation below will show). 

By Fact 3 and Fact 3''
the restriction of $\psi$ to $Z^{\#\#}$ is immersive outside
$E_0^{\#}$, and what we need to check is that this remains true on this 
divisor. The intersection $Z^{\#\#}\cap E_0^{\#}$ is a bunch
of projective lines -- the pre-images under the blow-up $\sigma$  
of the finite number of points inside $Z^{\#}\cap Ind(\varphi)$. 

Let $(p_1, H_1)$ be one of these points. We use the 
notations of the proof of Fact 6, with $n=3$. Let $F=0$ be the
equation of the surface $Z$, and $f(x_0,x_2,x_3)=0$ the equation we get 
by letting $x_1=1$. Our hypothesis on $H_1=T_{p_1}Z$ means that 
$\partial f/\partial x_0(0)=\partial f/\partial x_3(0)=0$, and we may    
suppose that $\partial f/\partial x_2(0)=1$. The local equations
of $Z^{\#}$ are 
$$f(x_0,x_2,x_3)=0, \quad  (\partial f/\partial x_2)y^0=
\partial f/\partial x_0, \quad  (\partial f/\partial x_2)y^3=
\partial f/\partial x_3.$$ 
If we write $f(x_0,x_2,x_3)=x_2+q(x_0,x_2,x_3)+$ {\it higher order terms},
where $q$ is quadratic, we deduce from these equations that the tangent 
space of $Z^{\#}$ at $(p_1, H_1)$ is given by:
$$\begin{array}{rcl}
dx_2 & = & 0, \\
dy^0 & = & q_{00}dx_0+q_{03}dx_3, \\
dy^3 & = & q_{03}dx_0+q_{33}dx_3.
\end{array}$$
In particular, the intersection with $Ind(\varphi)$ is transverse
if and only if $q_{33}\neq 0$. But in the system of local coordinates
on $Z$ given by $x_0$ and $x_3$, the local analytic equation of that 
hypersurface is, up to higher order terms, $x_2+q_{00}x_0^2
+2q_{03}x_0x_3+q_{33}x_3^2=0$, so that the condition $q_{33}\neq 0$
precisely means that the second fundamental form $II : (x_0,x_3)
\mapsto q_{00}x_0^2+2q_{03}x_0x_3+q_{33}x_3^2$ does not vanish 
identically on the line $x_0=0$. This is precisely condition $(3)$
at $p_1$. 

Now consider the strict transform of $Z^{\#}$, which means that
we let $x_0=ty^3$ and replace $x_0$ by $t$ in our system of local equations. 
The equations of the tangent space $T$ of $Z^{\#\#}$ at the point 
$(0,0,[1,t])$ over $(p_1, H_1)$ become:
$$\begin{array}{rcl}
dx_2 & = & 0, \\
dy^0 & = & tq_{00}dy^3+q_{03}dx_3, \\
dy^3 & = & tq_{03}dy^3+q_{33}dx_3.
\end{array}$$
But $\psi_*(dt,dx_2,dx_3,dy^0,dy^3)=(-tdx_2,dt,dx_3-tdy^0,dx_2,tdy^3)$. 
If $t\neq 0$, the kernel of $\psi_*$ is the line $dt=dx_2=dy^3=dx_3-tdy^0=0$, 
which is not contained in $T$. If $t=0$, this kernel is the plane 
$dt=dx_2=dx_3=0$, which again meets $T$ only at the origin. This 
implies our claim. \end{proof}

\bigskip
The series of conditions given by Propositions \ref{p17} and \ref{p18} look very restrictive. 
In particular, 
a ``generic'' smooth surface will have a curve (possibly reducible)
of bitangent planes, covering the whole of $\PP^3$. This will
not be compatible with condition $(2)$ of Proposition \ref{p17}, which 
is really a necessary condition for $\tilde{Z}$ to be smooth. 

Thus we   look for singular surfaces, but not 
too singular since they must remain quasi-smooth, at least outside
the hyperplane $H_0$. Let ${\check p}_0$ denote the set of hyperplanes passing through
$p_0$. 

\begin{theo} Let $Z$ be a Kummer quartic surface in $\PP^3$, in general position 
with respect to $p_0$ and $H_0$. Let $C=Z\cap H_0$ and $D=Z^*\cap {\check p}_0$, two general
hyperplane sections. The pull-backs of these curves to the K3 surface $Z^{\#}$ meet transversely
in twelve points. The surface $\tilde{Z}$ is isomorphic to $Z^{\#}$ blown up 
at these twelve points and is a smooth Legendrian surface in $\PP^5$. \end{theo}

Recall (see for example \cite {gh}) that a Kummer surface $Z$ in $\PP^3$ is a singular quartic 
surface in $\PP^3$ with exactly $16$ ordinary double points as 
singularities. In particular, it is quasi-smooth. But the 
property that will be most useful to us is that 
the dual of a Kummer surface also only has ordinary double points
as singularities. In fact the Kummer surface
is projectively isomorphic to its dual surface $Z^{\ast}$.  The 
double points of the dual surface define $16$ planes in $Z$ called
its {\em tropes}. These planes are (doubly) tangent to $Z$ along smooth 
conics, each of which contains exactly $6$ double points (this is the famous
$16_6$ configuration). 
The map $Z^{\#}\lra Z$ blows-up the $16$ double points of $Z$, 
and the dual map $Z^{\#}\lra Z^*$ contracts the pull-back of the 
tropes to the $16$ double points of $Z^*$. 

\begin{proof} We verify that the conditions of the previous two Propositions hold
for general $p_0$ and $H_0$. 
The only bitangent planes to $Z$ are the $16$ tropes, so the conditions of Proposition 
\ref{p17} clearly hold true in general. 
Since $Z$ has a finite number of singular points, which are ordinary double points, 
it is certainly quasi-smooth, which was condition (1) of Proposition \ref{p18}. 
If the hyperplane ${\check p}_0$ of ${\check\PP}^3$ does not contain any of the
sixteen singular points of $Z^*$, then the projection $Z^{\#}\lra Z^*$ is an isomorphism
above $Z^*\cap {\check p}_0$, and condition (2) is also verified. Finally, (3) is   again
a general position condition and will hold in general. We conclude that 
$Z^{\#\#}$ is smooth and isomorphic with $\tilde{Z}$.

Now the hyperplane section   $Z\cap H_0$, supposed to be   general, is a smooth quartic plane curve $C$ whose dual, by the Pl\"ucker
formulas,  is a curve 
of degree twelve. We conclude that the indeterminacy locus of 
$\varphi$ restricted to $Z$ is given by the twelve points on $C$ 
whose tangent line to $C$ hits $p_0$. The surface $Z^{\#\#}$ is therefore 
the $K3$-surface $Z^{\#}$, blown-up at the twelve corresponding points. 
\end{proof} 

We can give a precise meaning to the condition that $Z$ be in general
position with respect to $p_0$ and $H_0$. Namely, we need that:
\begin{enumerate}
\item $p_0$ is not contained in $Z$ nor in any of its tropes, 
\item $H_0$ is not tangent to $Z$ and contains none of its double points, 
\item none of the $28$  bitangents of the quartic curve $C=Z\cap H_0$ 
pass  through $p_0$, 
\item if $p\in C$ is such that the tangent line $T_pC$ contains $p_0$,
this line is not a bitangent, and is 
not contained in the kernel of the second fundamental 
form of $Z$ at $p$. 
\end{enumerate}  
These are all non-empty and open conditions. 

\medskip\noindent Note that since $c_2(Z^{\#})=24$, we get $c_2(\tilde{Z})=36$, 
and since $c_1(Z^{\#})=0$, $c_1(\tilde{Z})$ is minus the sum of the $12$ exceptional 
divisors $E_1,\ldots ,E_{12}$. In particular, $c_1^2=-12$ and $2ch_2=-84$. 

On the other hand, let $L$ (resp. $L'$) denote the pull-back to $Z^{\#}$ of the 
hyperplane divisor of $Z$ (resp. $Z^*$). By construction, the hyperplane 
class on $\tilde{Z}$ is $h = L+L'-E_1-\cdots -E_{12}$. We know that 
$L^2=(L')^2=4$, while $L.L'=12$, hence $h^2=4+2\times 12+4-12=20$
and $c_1.h=-12$. Finally, $2c_1h-3h^2=-2\times 12-3\times 20=-84$, 
in agreement with the identity $2ch_2=2c_1h-3h^2$ which we proved to hold for 
Legendrian surfaces.  
\medskip

\subsection{Relations with homalo\"{\i}dal polynomials}
In this final section we explore the relation of Bryant's method with the 
subadjoint varieties. 

Let $P$ be a homogeneous polynomial of degree $d$ in $n$ variables.
Denote by $Z_P\subset\PP^{n+1}$ the hypersurface of equation 
$x_0^{d-1}x_{n+1}=P(x_1,\ldots ,x_n)$.
After applying Bryant's birational map $\varphi$ to its tranform 
$Z_P^{\#}\in\FF_{1,n}(\CC^{n+1})$, we get a Legendrian, possibly singular 
variety $\tilde{Z}_P$. For $d\neq 2$,this variety can also be described as 
the image of the birational map 
$$\psi : [x_0,x]\in\PP^n\mapsto [x_0^d,  x_0^{d-1}x, x_0\partial P, P]\in\PP^{2n+1}.$$

\smallskip\noindent 
Again the problem is: when is $\tilde{Z}_P$ a smooth Legendrian variety ? 

\medskip \noindent {\it Example 1}. Let $P=x_1^3$, so that $Z_P$ is a cuspidal rational
plane cubic. Its dual $Z_P^*$ is again a cuspidal cubic, but $Z_P^{\#}$ is smooth, and 
$\tilde{Z}_P$ is a normal rational cubic curve in $\PP^3$.

\medskip \noindent {\it Example 2}. Let $P=x_1x_2...x_n$. One can check that $\tilde{Z}_P$
is singular for $n\ge 4$. For $n=3$, $\tilde{Z}_P=\PP^1\times\PP^1\times\PP^1$ is smooth. 

\medskip \noindent {\it Example 3}. Let $q$ be a quadratic form on $\CC^{n-1}$, 
and let $P=q(x_1,\ldots ,x_{n-1})x_n$. 
Then $\tilde{Z}_P\simeq\PP^1\times\QQ (q)$, where $\QQ(q)\subset\PP^n$ denotes the 
quadric hypersurface of equation $x_0x_n=q(x_1,\ldots ,x_{n-1})$. In particular, 
$\tilde{Z}_P$ is smooth if and only if $q$ is nondegenerate. 

\medskip \noindent We first note that the self-duality of $Z$ is a rather 
general phenomenon. This geneneralizes the well-known self-duality of smooth 
quadrics, which is the case where $d=2$. We refer to \cite{esb} for the terminology
in the next statement. 

\begin{prop} Let $P$ denote the unique (up to constant) relative invariant of minimal degree 
of an irreducible regular prehomogeneous space under some reductive Lie group. 
Then the hypersurface $Z$ is self dual. \end{prop}

\begin{proof} As noticed in \cite{esb}, the polynomial $P(\partial_1P,\ldots ,\partial_nP)$
is again a relative invariant, nonzero because of the regularity hypothesis. It must therefore 
be a nonzero multiple of $P^{d-1}$. In particular the homogeneous coordinates $
[y_0=(d-1)x_0^{d-2}x_{n+1},y_1=\partial_1P,\ldots, y_n=\partial_nP,y_{n+1}=x_0^d]$
are related by an  identity 
$$P(y_1,\ldots,y_n)=cP(x)^{d-1}=c(x_0^dx_{n+1})^{d-1}=(\frac{y_0}{d-1})^{d-1}y_{n+1}.$$
This proves the claim. \end{proof}

\medskip Let us try to understand when $\tilde{Z}_P$ can be 
smooth at the point $q_0=[0,\ldots ,0,1]$. We first note that its 
(reduced) tangent cone contains, for each $x$ such that $P(x)\neq 0$, the 
point $[0,0,\partial P(x),0]$. We can suppose that $x\mapsto\partial P(x)$ 
is linearly non-degenerate: otherwise, after a change of coordinates we can
suppose that $\partial P/\partial x_n=0$, hence $\tilde{Z}_P$ is linearly 
degenerate, hence a linear space. Then
the reduced tangent cone must contain the $n$-dimensional linear 
space $L=[0,0,\ast, \ast]$. 

Supposing $\tilde{Z}_P$ to be smooth, $L$ must co\"{\i}ncide with its tangent space
at $q$. Then the projection of $Z$ on $L$ with respect to the supplementary
space $[\ast,\ast,0,0]$ must be a local isomorphism in the complex topology. 
That is, the map $[x_0^d,  x_0^{d-1}x, x_0\partial P, P]\mapsto [x_0\partial P, P]$
must be a local isomorphism at $q$ -- in particular it must be injective. 
But suppose that we can find two noncolinear vector $x_1$ and $x_2$ 
with $P(x_1),P(x_2) \neq 0$, such that the vectors $\partial P(x_1)$ and 
$\partial P(x_2)$ are colinear. After multiplying them by a suitable constant, 
we may suppose that $P(x_1)^{-1}\partial P(x_1)=P(x_2)^{-1}\partial P(x_2)$. 
Then for $x_0$ small enough, $[x_0^d,  x_0^{d-1}x_1, x_0\partial P(x_1), P(x_1)]$
and $[x_0^d,  x_0^{d-1}x_2, x_0\partial P(x_2), P(x_2)]$ are two distinct points
in $\tilde{Z}$, arbitrarily close to $q$, but $[x_0\partial P(x_1), P(x_1)]$ and 
$[x_0\partial P(x_2), P(x_2)]$ coincides -- a contradiction with the smoothness
of $\tilde{Z}_P$ at $q$. We conclude that the rational map $[x]\in\PP^{n-1}
\mapsto [\partial P(x)]\in\PP^{n-1}$ must be injective on the open subset 
$P(x)\neq 0$. In particular, we have proved:

\begin{prop}
Let $Z_P\subset\PP^{n+1}$ be the hypersurface $x_0^{d-1}x_{n+1}=P$, and 
suppose that $\tilde{Z}_P$ is a smooth Legendrian variety. Then $P$ must be 
homalo\"{\i}dal.
\end{prop}

\begin{coro} If $d=3$ and $\tilde{Z}_P$ is a smooth Legendrian variety, then $P$ must be 
the determinant of a semisimple Jordan algebra of rank three, and $\tilde{Z}_P$ is a 
subadjoint variety. \end{coro}

\begin{proof} That a homalo\"{\i}dal polynomial of degree three must be the determinant of a semisimple 
Jordan algebra of rank three is due to Etingof, Kazhdan and Polischuk \cite{ekp}, see also 
\cite{chaput}. The fact that the resulting varieties $\tilde{Z}_P$ are exactly the subadjoint
varieties follows from \cite{LMclass}, where they were constructed precisely that way. \end{proof}

\end{document}